\documentclass[a4paper,11pt]{article}
\usepackage[hidelinks]{hyperref}
\hypersetup{
	colorlinks=true,   
	linktoc=all,        
	linkcolor=black,     
	citecolor=blue,}     

\usepackage{tikz}
\usepackage{graphicx}
\usepackage{subfigure}
\usepackage{amssymb}
\usepackage{latexsym,bm}
\usepackage{multicol}
\usepackage{indentfirst}
\usepackage{amssymb,amsfonts}
\usepackage{amsmath}
\usepackage{setspace}
\newcommand{\ignore}[1]{}

\newtheorem{theorem}{Theorem}[section]
\newtheorem{lemma}[theorem]{Lemma}

\newtheorem{conjecture}[theorem]{Conjecture}

\begin{document}
	\date{}
	\begin{spacing}{1.03}
		\title{
Phase transitions of the Erd\H{o}s-Gy\'{a}rf\'{a}s function}
		
		\author{Xinyu Hu,\footnote{Data Science Institute, Shandong University, Jinan, 250100, P.~R.~China. Email: {\tt huxinyu@sdu.edu.cn}. }
			\;\; \; Qizhong Lin,\footnote{Center for Discrete Mathematics, Fuzhou University, Fuzhou, 350108, P.~R.~China. Email: {\tt linqizhong@fzu.edu.cn}. Supported in part  by National Key R\&D Program of China (Grant No. 2023YFA1010202), NSFC (No.\ 12171088) and NSFFJ (No. 2022J02018).}
\;\; \; Xin Lu,\footnote{
Zhengdong New District Yongfeng School, Zhengzhou, 451460, P. R. China. Email: {\tt lx20190902@163.com}.}
\;\; \;Guanghui Wang \footnote{
School of Mathematics, Shandong University, Jinan, 250100, P. R. China. Email: {\tt ghwang@sdu.edu.cn}. Supported in part  by State Key Laboratory of Cryptography and Digital Economy Security.} 
}
		\maketitle

		\begin{abstract}
Given positive integers $p,q$. For any integer $k\ge2$, an edge coloring of the complete $k$-graph $K_n^{(k)}$ is said to be a $(p,q)$-coloring if every copy of $K_p^{(k)}$ receives at least $q$ colors. The Erd\H{o}s-Gy\'{a}rf\'{a}s function $f_k(n,p,q)$ is the minimum number of colors  that are needed for $K_n^{(k)}$ to have a $(p,q)$-coloring.

Conlon, Fox, Lee and Sudakov (\emph{IMRN, 2015}) conjectured that for any positive integers $p, k$ and $i$  with $k\ge3$ and $1\le i<k$, $f_k(n,p,{{p-i}\choose{k-i}})=(\log_{(i-1)}n)^{o(1)}$, where $\log_{(i)}n$ is an iterated $i$-fold logarithm in $n$. It has been verified to be true for $k=3, p=4, i=1$ by Conlon et. al (\emph{IMRN, 2015}), for $k=3, p=5, i=2$ by Mubayi (\emph{JGT, 2016}), and for all $k\ge 4, p=k+1,i=1$ by B. Janzer and O. Janzer (\emph{JCTB, 2024}).
In this paper, we give new constructions and show that this conjecture holds for infinitely many new cases, i.e., it holds for all $k\ge4$, $p=k+2$ and $i=k-1$.


\end{abstract}

\section{Introduction}
A $k$-uniform hypergraph $H$ ($k$-graph for short) with vertex set $V(H)$ is a collection of $k$-element subsets of $V(H)$. We write $K_n^{(k)}$
for the complete $k$-graph on an $n$-element vertex set. Ramsey theorem \cite{R} implies that for any integers $n_1,\dots,n_q$, there exists a minimum integer, now called Ramsey number $N=r_k(n_1,\dots,n_q)$, such that any $q$-coloring of edges of the complete $k$-graph $K_N^{(k)}$ contains a $K_{n_i}^{(k)}$ in the $i$th color for some $i\in[q]$. We will use the simpler notation $r_k(n;q)$ if $n_i=n$ for all $i$.

A $(p,q)$-coloring of $K_n^{(k)}$ is an edge-coloring of $K_n^{(k)}$ that gives every copy of $K_p^{(k)}$
at least $q$ colors. Let $f_k(n,p,q)$ be the minimum number of colors in a $(p,q)$-coloring of $K_n^{(k)}$. The function $f_k(n,p,q)$ can be seen as a generalization of the usual Ramsey function. Indeed, when $q=2$, we know that
\begin{align}\label{relat}
\text{$f_k(n,p,2)=\ell$ if and only if  $r_k(p;\ell)> n$ and $r_k(p;\ell-1)\le n.$}
\end{align}
 Therefore, when determining $f_k(n,p,q)$, we are generally interested in $q\ge 3$. For simplicity, we write $f(n,p,q)$ when $k=2$.

Erd\H{o}s and Shelah \cite{R-S-1,R-S-2} initiated to determine $f(n,p,q)$ for fixed $p$ and $q$ where $2\le q\le{\binom p 2}$. Subsequently, Erd\H{o}s and Gy\'{a}rf\'{a}s \cite{E-G-1} systematically studied this function. Since the function $f(n,p,q)$ is increasing in $q$, we are interested in the transitions of $f(n,p,q)$ as $q$ increases. It is clear that $f(n,p,2)\le f(n,3,2)=O(\log n)$ by noting $r_2(3;t) > 2^t$, while $f(n,p,{\binom p 2})={\binom n 2}$ for $p\ge4$. In particular, Erd\H{o}s and Gy\'{a}rf\'{a}s \cite{E-G-1} proved that for $p\ge3$,
\begin{align}\label{low-pp}
n^{1/(p-2)}-1\le f(n,p,p) \le O(n^{2/(p-1)}),
\end{align} which implies that $f(n,p,q)$ is polynomial in $n$ for any $q \ge p$. 
Erd\H{o}s and Gy\'{a}rf\'{a}s asked whether
$f(n,p,p-1) = n^{o(1)}$ for all fixed $p\ge4$. If the answer is yes, then $p-1$ is the maximum $q$ such that $f(n,p,q)$ is subpolynomial in $n$ by noting $f(n,p,p)=\Omega(n^{1/(p-2)})$.
The first case was verified by Mubayi \cite{M-1} from an elegant construction, indeed, Mubayi established
$$f(n,4,3)=e^{O(\sqrt{\log n})}.$$
 The best lower bound $f(n,4,3)=\Omega(\log n)$ is due to Fox and Sudakov \cite{F-S-1}, improving that by Kostochka and Mubayi \cite{K-M-1}. Applying the same construction of \cite{M-1}, Eichhorn and Mubayi \cite{E-M-1} obtained that $f(n,5,4)=e^{O(\sqrt{\log n})}$.
Finally, Conlon, Fox, Lee and Sudakov \cite{C-F-L-S-1} answered the question in the affirmative. In fact, they showed that for any fixed $p\ge4$,
\[
f(n,p,p-1)\le e^{(\log n)^{1-1/(p-2)+o(1)}} = n^{o(1)}.
\]
Moreover, the exponent $1/(p-2)$ in the lower bound (\ref{low-pp}) was shown
to be sharp for $p=4$ by Mubayi \cite{M-2} and also for $p=5$ by Cameron and Heath \cite{C-H-1} via explicit constructions. Recently, Cameron and Heath \cite{C-H-2} showed that $f(n,6,6)\le n^{1/3+o(1)}$ and $f(n,8,8)\le n^{1/4+o(1)}$.

The first nontrivial hypergraph case is $f_3(n,4,3)$, which has tight connections to Shelah's breakthrough proof \cite{S-1} of primitive recursive bounds for the Hales-Jewett numbers. Answering a question of Graham, Rothschild and Spencer \cite{G-R-S-1}, Conlon, Fox, Lee and Sudakov \cite{C-F-L-S-2} showed that $$f_3(n,4,3)=n^{o(1)}.$$ In general, using a variant of the pigeonhole argument for hypergraph Ramsey
numbers due to Erd\H{o}s and Rado, Conlon et. al \cite{C-F-L-S-2} proved that for any fixed positive integers $p,k,i$, there exists a constant $c>0$ such that
 $$f_k\left(n,p,{{p-i}\choose{k-i}}+1\right)=\Omega(\log_{(i-1)}n^{c}),$$where $\log_{(0)}(x) = x$ and $\log_{(i+1)} x = \log \log_{(i)} x$ for $i\ge0$. They \cite[Problem 6.2]{C-F-L-S-2} (from the perspective of the inverse problem) further proposed a variety of basic questions about the Erd\H{o}s-Gy\'{a}rf\'{a}s function $f_k(n,p,q)$ as follows.
\begin{conjecture}[Conlon, Fox, Lee and Sudakov \cite{C-F-L-S-2}]\label{conj}
Let $p, k$ and $i$ be positive integers with $k\ge3$ and $1\le i<k$, $$f_k\left(n,p,{{p-i}\choose{k-i}}\right)=(\log_{(i-1)}n)^{o(1)}.$$

\end{conjecture}

We can see that when $k=2$, it is precisely the Erd\H{o}s-Gy\'{a}rf\'{a}s problem. For the case where $k=3, p=4, i=1$, we know  \cite{C-F-L-S-2} that Conjecture \ref{conj} holds. The next case $k=3, p=5, i=2$
 was verified by Mubayi \cite{M-3}, who indeed showed that
$f_3(n,5,3)=e^{O(\sqrt{\log \log n})}=(\log n )^{o(1)}$.
Recently, B. Janzer and O. Janzer \cite{J-J-1} showed that $f_k(n,k+1,k)=n^{o(1)}$ for all $k\ge 4$, together with that obtained in \cite{C-F-L-S-2} implying that Conjecture \ref{conj} holds for all $k\ge 3$, $p=k+1$ and $i=1$. 
 We refer the reader to \cite{C-F-S,M-S-1} for two nice surveys on this topic.




In this paper, we obtain an upper bound of $f_k(n,k+2,3)$ for all $k\ge 4$ as follows. Together with the case where $k=3$ due to Mubayi \cite{M-3}, we know that Conjecture \ref{conj} holds for all $k\ge 3$, $p=k+2$ and $i=k-1$.
\begin{theorem}\label{center-2}
For any fixed integer $k\ge 4$, $f_k(n,k+2,3)=e^{O(\sqrt{\log_{(k-1)} n})}=(\log_{(k-2)} n)^{o(1)}.$

\end{theorem}

Our construction is based on the Mubayi's coloring in \cite{M-3}, and we define the auxiliary color mapping to construct a $(k+2,3)$-coloring. Moreover, we using the stepping up technique of Erd\H{o}s and Hajnal.

The organization of the paper is as follows. In section 2 we will give some notation and basic properties. In section 3 we will give the coloring constructions. More precisely, in subsection 3.1 we will recall the explicit edge-colorings constructed by Mubayi. In subsection 3.2 we will prove Theorem \ref{center-2} for the case $k=4$, i.e., $f_4(n,6,3)=e^{O(\sqrt{\log\log\log n})}=(\log\log n)^{o(1)}.$ In subsection 3.3 we will show Theorem \ref{center-2}.



\section{Notation and basic properties}
In this paper, we will apply several variants of the Erd\H{o}s-Hajnal stepping-up lemma. Given some integer number $N$, let $V=\{0,1\}^N$. The vertices of $V$ are naturally ordered by the integer they represent in binary, so for $a,b\in V$ where $a=(a(1),\ldots,a(N))$ and $b=(b(1),\ldots,b(N))$, $a<b$ \textbf{iff} there is an $i$ such that $a(i)=0$, $b(i)=1$, and $a(j)=b(j)$ for all $1\le j<i$. In other words, $i$ is the first position (minimum index) in which $a$ and $b$ differ. For $a \not = b$, let $\delta(a,b)$ denote the minimum $i$ for which $a(i) \not = b(i)$. Given any vertices subset $S=\{a_1,\ldots,a_{r}\}$ of $V$ with $a_1<\cdots<a_r$, we always write for $1\le s<t\le r$, $$\delta_{st}=\delta(a_s,a_{t}).$$ If $t=s+1$, we will use the simpler notation $\delta_s=\delta(a_s,a_{s+1})$. We say that $\delta_s$ is a \emph{local minimum} if $\delta_{s-1}>\delta_{s}<\delta_{s+1}$, a \emph{local maximum} if $\delta_{s-1}<\delta_{s}>\delta_{s+1}$, and a \emph{local extremum} if it is either a local minimum or a local maximum. For convenience, we write $\delta(S)=\{\delta_s\}_{s=1}^{r-1}$.


\medskip

We have the following stepping-up properties, see in \cite{G-R-S-1}.

\begin{description}

\item[Property A:] For every triple $a < b < c$, $\delta(a,b) \not = \delta(b,c)$ .

\item[Property B:] For $a_1 <a_2< \cdots < a_r$, $\delta_{1r}=\delta(a_1,a_{r}) = \min_{1 \leq i \leq r-1}\delta_i$.

\end{description}

Since $\delta_{s-1}\neq\delta_{s}$ for every $s$, every nonmonotone sequence $\{\delta_s\}_{s=1}^{r-1}$ has a local extremum.
\medskip

We will also use the following stepping-up properties, which are easy consequences of Properties A and B, see e.g. \cite{F-H-L-L}, and we include the proofs for completeness.

\begin{description}
\item[Property C:]  For $\delta_{1r}=\delta(a_1,a_{r}) = \min_{1 \leq i \leq r-1}\delta_i$, there is a unique $\delta_i$ which achieves the minimum.
 \end{description}

\noindent
{\bf Proof.} Suppose for some $s<t$,
$\delta(a_s, a_{s+1})=\delta(a_{t}, a_{t+1})=\min_{1 \leq i \leq r-1} \delta_i.$
Then, by Property B, $\delta(a_s, a_{t})=\delta(a_{t}, a_{t+1})$, contradicting Property A.
\hfill$\Box$

\begin{description}
\item[Property D:] For every 4-tuple $a_1 <a_2<a_3 < a_4$, if $\delta_1 < \delta_2$, then $\delta_1 \neq \delta_3$.
 \end{description}
\noindent
{\bf Proof.} Otherwise, suppose $\delta_1 = \delta_3$.
Then, by Property B, $\delta(a_1, a_{3})=\delta_1 =\delta_3=\delta(a_{3}, a_{4})$. This contradicts Property A since $a_{1}<a_{3}<a_{4}$.
\hfill$\Box$

\section{The coloring constructions}\label{pre}

In this section, we will prove our main results through constructing the suitable $(p,q)$-colorings, see subsections \ref{con-6-3} and \ref{con-gen}. For clarity, we write $\chi_i$ ($i\ge2$) for the edge-coloring of the complete $i$-graph $H_i$ on $N_i$ (ordered) vertices, where
\[
N_2={m\choose t}, \;\;\text{and} \;\; N_{i+1}=2^{N_i} \;\;\text{for} \;\; i\ge2.
\]

\subsection{Mubayi's-colorings}
We first recall the explicit edge-coloring $\chi_2$ constructed by Mubayi \cite{M-1}, from which we know that $f(n,4,3)=e^{O(\sqrt{\log n})}$.

\medskip
\noindent
\textbf{Construction of $\chi_2$:} Given integers $t<m$ and $N_2={{m}\choose {t}}$, let $V(K_{N_2} )$ be the set of $0/1$ vectors of length $m$ with exactly $t$ 1's. Write $v=(v(1), \ldots,v(m))$ for a vertex of $K_{N_2}$. The vertices
are naturally ordered by the integer they represent in binary, so $v<w$ iff $v(i) = 0$ and $w(i) = 1$ where $i$ is the first position (minimum integer) in which $v$ and $w$ differ. By considering vertices as characteristic vectors of sets, we may assume $V(K_{N_2} )={{[m]}\choose {t}}$. For each $B\in {{[m]}\choose {t}}$, let $f_B:2^B \rightarrow[2^t]$ be a bijection. Given vectors $v<w$
that are characteristic vectors of sets $S<T$, let
\begin{align*}
   & c_1(vw)=\min\{i:v(i)=0, w(i)=1\}, \\
   & c_2(vw)=\min\{j:j>c_1(vw), v(j)=1, w(j)=0\}, \\
   & c_3(vw)=f_S(S \cap T ), \\
   & c_4(vw)=f_T (S\cap T ).
\end{align*}
Finally, define
\[
\chi_2 (vw)=(c_1(vw), c_2(vw), c_3(vw), c_4(vw)).
\]
If ${N_2}$ is not of the form ${{m}\choose {t}}$, then let ${N_2}'\ge {N_2}$ be the smallest integer of this form, color ${{[{N_2}']}\choose {2}}$ as described above, and restrict the coloring to ${{[{N_2}]}\choose {2}}$.
It is known \cite{M-1,M-2} that $\chi_2$ is both a $(3, 2)$ and $(4, 3)$-coloring of $K_{N_2}$ (We only need the
first and fourth coordinates of $\chi_2$ for this) and, for suitable choice of $m$ and $t$, it uses $e^{O(\sqrt{\log {N_2}})}$ colors for all ${N_2}$.
Therefore, we have \cite{M-1} that $f(n,4,3)=e^{O(\sqrt{\log n})}$.

\medskip

Now we recall the edge-coloring $\chi_3$ due to Mubayi \cite{M-3}.

\medskip\noindent
\textbf{Construction of $\chi_3$:} Given a copy of $K_{N_2}$ on $[{N_2}]$ and the edge-coloring $\chi_2$, and let $N_3=2^{N_2}$. We produce an edge-coloring $\chi_3$ of $H_3$ on $\{0,1\}^{N_2}$ as follows. Order the vertices of $H_3$ according to the integer that they represent in binary. For an edge $(a_i,a_j,a_k)$ with $a_i<a_j<a_k$, then $\delta_{ij}\neq \delta_{jk}$ from Property A.
Let
\begin{align*}
\chi_3(a_i,a_j,a_k)=(\chi_2(\delta_{ij},\delta_{jk}),\delta_{ijk}),
\end{align*}
 where $\delta_{ijk}$ equals 1 if $\delta_{ij} <\delta_{jk}$ and $-1$ otherwise. Since $\chi_2$ is an edge-coloring of $K_{N_2}$ with $e^{O(\sqrt{\log {N_2}})}$ colors, we obtain that $\chi_3$ is an edge-coloring of $H_3$ with $e^{O(\sqrt{\log \log {N_3}})}$ colors as promised. 

From the following property of $\chi_3$, we know that $f_3(n,5,3)=e^{O(\sqrt{\log \log n})}$.
\begin{lemma}[Mubayi \cite{M-3}]\label{(5,3)-coloring}
$\chi_3$ is a $(5,3)$-coloring of $H_3$.
\end{lemma}

We also need the following property of $\chi_3$.
\begin{lemma}\label{(4,2)-coloring}
$\chi_3$ is a $(4,2)$-coloring of $H_3$.
\end{lemma}

\noindent
\textbf{Proof.} Suppose, for contradiction that $Y_3=\{a_1,\ldots,a_4\}$ with $a_1<a_2<a_3<a_4$ are four vertices of $H_3$ forming a monochromatic $K_{4}^{(3)}$. Recall that $\delta_i=\delta(a_i,a_{i+1})$ for $i\in[3]$. Let $\delta=\min_{1\le i\le 3} \delta_i$. We know $\delta_i\neq \delta_{i+1}$ for $i\in[2]$ from Property A. Suppose first that $\delta=\delta_1$, then $\delta_1<\delta_2$. If $\delta_2<\delta_3$, then the $K_3$ on $\{\delta_1,\delta_2,\delta_3\}$ has two colors since $\chi_2$ is a $(3,2)$-coloring and this gives two colors to the edges of $H_3$ within
$\{a_1,\ldots,a_4\}$. Thus, we assume $\delta_2>\delta_3$. Note that $\delta_{123}=1$ as $\delta_1<\delta_2$, and $\delta_{234}=-1$ as $\delta_2>\delta_3$, then $\chi_3(a_1,a_2,a_3)\neq\chi_3(a_2,a_3,a_4)$, and so the $K_{4}^{(3)}$ on $\{a_1,\ldots,a_4\}$ have two colors. The case for $\delta=\delta_3$ is similar. Now suppose $\delta=\delta_2$, then $\delta_1>\delta_2<\delta_3$. Note that $\delta_{123}=-1$ as $\delta_1>\delta_2$, and $\delta_{234}=1$ as $\delta_2<\delta_3$, then $\chi_3(a_1,a_2,a_3)\neq\chi_3(a_2,a_3,a_4)$, and so the $K_{4}^{(3)}$ on $\{a_1,\ldots,a_4\}$ have two colors. \hfill$\Box$

\subsection{Construction of $(6,3)$-coloring}\label{con-6-3}
In this subsection, we aim to construct a $(6,3)$-coloring $\chi_4$ of $K_{N_4}^{(4)}$ on $N_4$ vertices from the coloring $\chi_3$ defined in the last subsection, from which we will show $f_4(n,6,3)=e^{O(\sqrt{\log\log\log n})}$.

Given a copy of $H_3$ on $[N_3]$ and the edge-coloring $\chi_3$, we will produce an edge-coloring $\chi_4$ of the complete $4$-graph $H_4:=K_{N_4}^{(4)}$ on $\{0,1\}^{N_3}$ as follows.  Order the vertices of $H_4$ according
to the integer that they represent in binary. Given an edge $e=(a_i,a_j,a_k,a_\ell)$ with $a_i<a_j<a_k<a_\ell$, we first define one auxiliary color mapping $\varphi_3$. Recall that $\delta_{st}=\delta(a_s,a_t)$ for $a_s<a_t$, and $\delta_{ij}\neq\delta_{k\ell}$ if $\delta_{ij}<\delta_{jk}>\delta_{k\ell}$ from Property D. Define $\delta(e)=\{\delta_{ij},\delta_{jk},\delta_{k\ell}\}$.
Let

\[
\varphi_3(\delta_{ij},\delta_{jk},\delta_{k\ell}) =
\begin{cases}
(I,0), & \text{if } \delta_{ij}<\delta_{jk}<\delta_{k\ell}, \\
(D,0), & \text{if } \delta_{ij}>\delta_{jk}>\delta_{k\ell},\\
(A,0), & \text{if } \delta_{ij}>\delta_{jk}<\delta_{k\ell},\\
(B,+), & \text{if } \delta_{ij}<\delta_{jk}>\delta_{k\ell} ~~\text {and} ~~\delta_{ij}<\delta_{k\ell},\\
(B,-), & \text{if } \delta_{ij}<\delta_{jk}>\delta_{k\ell}~~\text {and}~ ~\delta_{ij}>\delta_{k\ell},
\end{cases}
\]where $I,D,A,B$ have no inherent meaning and are only used to distinguish colors.

Now, we define $$\chi_4((a_i,a_j,a_k,a_\ell))=(\chi_3(\delta_{ij},\delta_{jk},\delta_{k\ell}),\varphi_3(\delta_{ij},\delta_{jk},\delta_{k\ell})).$$

Recall that $\chi_3$ is an edge-coloring of $K_{N_3}^{(3)}$ with $e^{O(\sqrt{\log \log {N_3}})}$ colors and ${N_4}=2^{N_3}$, then $\chi_4$ is an edge-coloring of $K_{N_4}^{(4)}$ with $5e^{O(\sqrt{\log \log {N_3}})}=e^{O(\sqrt{\log \log \log {N_4}})}$ colors as desired. Moreover, extending this construction to all ${N_4}$ is trivial by considering the smallest ${N_4}'\ge {N_4}$ which is a power of $2$, coloring ${{[{N_4}']}\choose{4}}$
and restricting to ${{[N_4]}\choose{4}}$.

\begin{lemma}\label{(6,3)-coloring}
Let $\chi_4$ be the edge-coloring of $H_4$ defined as above. Then $\chi_4$ is a $(6,3)$-coloring.
\end{lemma}
\textbf{Proof.} Let $X_4=\{a_1,\ldots ,a_6\}$ be vertices of $H_4$ forming a $K_6^{(4)}$ with $a_1<\cdots <a_6$, we will show that there has at least three colors in edges of $X_4$. Recall that $\delta_i=\delta(a_i,a_{i+1})$ for $i\in[5]$, and $\delta_{ij}=\delta(a_i,a_j)$ for $1\le i<j\le 6$. Let $\delta=\min_{1\le i\le 5} \delta_i$.  It follows from Property A and Property C that this minimal is uniquely achieved, and $\delta_i\neq \delta_{i+1}$ for $i\in[5]$. Let $e_1=(a_1,\ldots,a_4)$, $e_2=(a_2,\ldots,a_5)$, and $e_3=(a_3,\ldots,a_6)$ be the three edges in $X_4$, then $\delta(e_1)=\{\delta_1,\delta_2,\delta_3\}$, $\delta(e_2)=\{\delta_2,\delta_3,\delta_4\}$, and $\delta(e_3)=\{\delta_3,\delta_4,\delta_5\}$.

Suppose that $\chi_4(e_1)$, $\chi_4(e_2)$, and $\chi_4(e_3)$ are distinct to each other, then this gives three colors to the edges in $X_4$ and we need do nothing. Therefore, there are at least two of them that are equal. We split the proof into two cases as follows.

\medskip

{\em Case 1:} $\chi_4(e_1)=\chi_4(e_2)$, or $\chi_4(e_2)=\chi_4(e_3).$

\medskip
Suppose first that $\chi_4(e_1)=\chi_4(e_2)$, then $\varphi_3(\delta_1,\delta_2,\delta_3)=\varphi_3(\delta_2,\delta_3,\delta_4)$, which implies that $\{\delta_1,\delta_2,\delta_3,\delta_4\}$ is monotone. Without loss of generality, we may assume that $\delta_1<\cdots<\delta_4$. If $\delta_4<\delta_5$, then the $K_5^{(3)}$ on $\{\delta_1,\ldots,\delta_5\}$ has three colors since $\chi_3$ is a $(5,3)$-coloring and this gives at least three colors to the edges in $X_4$ as desired. Thus, $\delta_4>\delta_5$ from Property A, then the $K_4^{(3)}$ on $\{\delta_1,\ldots,\delta_4\}$ has two colors since $\chi_3$ is a $(4,2)$-coloring from Lemma \ref{(4,2)-coloring} and this gives two colors to the edges of $X_4$ within $\{a_1,\ldots,a_5\}$ and the $\varphi_3$-coordinate are $(I,0)$. Moreover, $\varphi_3(\delta_3,\delta_4,\delta_5)\in\{(B,+),(B,-)\}$, so we again have at least three colors on $X_4$. In the second case $\chi_4(e_2)=\chi_4(e_3)$, similar as above, we have at least three colors on $X_4$.

\medskip

\emph{Case 2:} $\chi_4(e_1)=\chi_4(e_3).$

\medskip
For this case, we have $\varphi_3(\delta_1,\delta_2,\delta_3)=\varphi_3(\delta_3,\delta_4,\delta_5)$.
Suppose that $\{\delta_1,\delta_2,\delta_3\}$ is monotone, and we may assume that it is monotone increasing without loss of generality. Then, $\{\delta_3,\delta_4,\delta_5\}$ is also monotone increasing. Therefore, the $K_5^{(3)}$ on $\{\delta_1,\ldots,\delta_5\}$ has three colors since $\chi_3$ is a $(5,3)$-coloring and this gives at least three colors to the edges in $X_4$ as desired.
So we may assume that $\{\delta_1,\delta_2,\delta_3\}$ is not monotone. From Property A, there are two cases.

If $\delta_1>\delta_2<\delta_3$, then $\delta_3>\delta_4<\delta_5$ since $\varphi_3(\delta_1,\delta_2,\delta_3)=\varphi_3(\delta_3,\delta_4,\delta_5)$, implying
 $$\delta_1>\delta_2<\delta_3>\delta_4<\delta_5.$$
From Property D, $\delta_2\neq \delta_4$. If $\delta_2> \delta_4$, then set $e_1'=(a_1,a_2,a_4,a_5)$, and $\delta(e_1')=\{\delta_1,\delta_2,\delta_4\}$ since $\delta_{24}=\min\{\delta_2,\delta_3\}=\delta_2$ by noting Property C, and so the $\varphi_3$-coordinate of $\chi_4(e_1')$ is $(D,0)$. Note that the $\varphi_3$-coordinate of $\chi_4(e_1)$ is $(A,0)$, and the $\varphi_3$-coordinate of $\chi_4(e_2)$ is $(B,-)$. Thus, we have at least three colors on $X_4$. Therefore, we assume $\delta_2< \delta_4$ and define $e_1'=(a_2,a_4,a_5,a_6)$, by a similar argument as above, we have at least three colors on $X_4$ as desired.

If $\delta_1<\delta_2>\delta_3$, then $\delta_3<\delta_4>\delta_5$, and so
$$\delta_1<\delta_2>\delta_3<\delta_4>\delta_5.$$
 From Property D, $\delta_1\neq \delta_3$. We may assume that $\delta_1>\delta_3$ without loss of generality, then $\delta_3>\delta_5$ since $\varphi_3(\delta_1,\delta_2,\delta_3)=\varphi_3(\delta_3,\delta_4,\delta_5)$. Let $e_3'=(a_2,a_3,a_4,a_6)$, then $\delta(e_3')=\{\delta_2,\delta_3,\delta_5\}$ since $\delta_{46}=\min\{\delta_4,\delta_5\}=\delta_5$ by noting Property C, and so the $\varphi_3$-coordinate of $\chi_4(e_3')$ is $(D,0)$. Note that the $\varphi_3$-coordinate of $\chi_4(e_1)$ is $(B,-)$, and the $\varphi_3$-coordinate of $\chi_4(e_2)$ is $(A,0)$. Thus, we have at least three colors on $X_4$.

This completes the proof.\hfill$\Box$


\subsection{Construction of $(k+2,3)$-coloring}\label{con-gen}


We will use induction on $k\ge 4$ to show Theorem \ref{center-2} via the suitable $(p,q)$-colorings. The based case is just proving in Subsection 3.2. For the inductive step, we assume that Theorem \ref{center-2} holds for $k-1$ with $k\ge5$, i.e., $f_{k-1}(n,k+1,3)=e^{O(\sqrt{\log_{(k-2)} n})}$. We aim to show $f_{k}(n,k+2,3)=e^{O(\sqrt{\log_{(k-1)} n})}$.

Given a copy of $K_{N_{k-1}}^{(k-1)}$ on $[N_{k-1}]$ and the edge-coloring $\chi_{k-1}$ from the induction hypothesis, we will produce an edge-coloring $\chi_k$ of the complete $k$-graph $H_k:=K_{N_{k}}^{(k)}$ on $\{0,1\}^{N_{k-1}}$ as follows.  Order the vertices of $H_k$ according
to the integer that they represent in binary. Given an edge $e=(a_{i_1},\ldots,a_{i_k})$ with $a_{i_1}<\cdots<a_{i_k}$, we
first define auxiliary color mapping $\varphi_{k-1}$. Recall that $\delta_{st}=\delta(a_s,a_t)$ for all $a_s<a_t$, and
$\delta(e)=\{\delta_{i_1i_2},\ldots,\delta_{i_{k-1}i_{k}}\}$. 
Note that every nonmonotone sequence has a local extremum, let $\delta_{\Delta}$ be the \textbf{first local extremum} of $\delta(e)$ if $\delta(e)$ is nonmonotone, where $\Delta\in \{i_2i_3,\ldots, i_{k-2}i_{k-1}\}$.


Let
\[
\varphi_{k-1}(\delta_{i_1i_2},\ldots,\delta_{i_{k-1}i_{k}}) =
\begin{cases}
(I,0), & \text{if } \delta(e) ~\text {is monotone increase,} \\
(D,0), & \text{if } \delta(e) ~\text {is monotone decrease.}\\
& \text{Otherwise, let $\ell\in[2,k-2]$ be the minimum index}
\\
& \text{such that $\delta_{\Delta}:=\delta_{i_{\ell}i_{\ell+1}}$ is a local extremum,}
\\
(A_\ell,0), & \text{if $\delta_{\Delta}$ is a local minimum},\\
(B_\ell,+), & \text{if $\delta_{\Delta}$ is a local maximum and $\delta_{i_{\ell-1}i_\ell}<\delta_{i_{\ell+1}i_{\ell+2}}$},\\
(B_\ell,-), & \text{if $\delta_{\Delta}$ is a local maximum and $\delta_{i_{\ell-1}i_\ell}>\delta_{i_{\ell+1}i_{\ell+2}}$}.\\
\end{cases}
\]
Now, we define $$\chi_k(a_{i_1},\ldots,a_{i_k})=(\chi_{k-1}(\delta_{i_1i_2},\ldots,\delta_{i_{k-1}i_k}),\varphi_{k-1}(\delta_{i_1i_2},\ldots,\delta_{i_{k-1}i_k})).$$
Since $\chi_{k-1}$ is an edge-coloring of $K_{N_{k-1}}^{(k-1)}$ with $e^{O(\sqrt{\log_{(k-2)}{N_{k-1}}})}$ colors from the induction hypothesis and ${N_k}=2^{N_{k-1}}$, $\chi_k$ is an edge-coloring of $K_{N_k}^{(k)}$ with $$(3k-7)e^{O(\sqrt{\log_{(k-2)}{N_{k-1}}})}=e^{O(\sqrt{\log_{(k-1)} {N_k}})}$$ colors as desired. Moreover, extending this construction to all ${N_k}$ is trivial by considering the smallest ${N_k}'\ge {N_k}$ which is a power of $2$, coloring ${{[{N_k}']}\choose{k}}$
and restricting to ${{[N_k]}\choose{k}}$.

Lemma \ref{(4,2)-coloring} tells us that $\chi_3$ is also a $(4,2)$-coloring. For general $k'\in [3,k]$, we can inductively prove the following lemma. 
\begin{lemma}\label{(k+1,2)-coloring}
For any $k\ge3$, $\chi_{k}$ is a $(k+1,2)$-coloring of $H_{k}$.
\end{lemma}
\textbf{Proof.} The proof is by using the induction on $k\ge3$. The base case where $k=3$ holds from Lemma \ref{(4,2)-coloring}. In general, for $k>3$, suppose the assertion holds for $k-1$ and we will show that it also holds for $k$.

Suppose to the contrary that $Y_{k}=\{a_1,a_2,\ldots,a_{k+1}\}$ forms a monochromatic $K_{k+1}^{(k)}$ in $H_{k}$. Consider two edges $e_1=(a_1,\ldots,a_{k})$ and $e_2=(a_2,\ldots,a_{k+1})$. Recall that $\delta_i=\delta(a_i,a_{i+1})$ for $i\in[k]$, $\delta_{ij}=\delta(a_i,a_j)$ for $1\le i<j\le k+1$, and $\delta(Y_k)=\{\delta_1,\dots,\delta_k\}$. If $\delta(Y_k)$ is monotone, then the $K_{k}^{(k-1)}$ on $\delta(Y_k)$ has two colors by noting $\chi_{k-1}$ is a $(k,2)$-coloring from the induction hypothesis and gives two colors to the edges of $H_{k}$ within $Y_{k}$. Thus, $\delta(Y_k)$ is nonmonotone, and so let $\delta_\ell$ be the first local extremum, where $\ell\in [2,k-1]$. We may assume that $\delta_\ell$ is a local maximum without loss of generality, i.e.,
$$\delta_1<\cdots<\delta_{\ell-1}<\delta_\ell>\delta_{\ell+1}\cdots.$$

If $\ell\in [3,k-1]$, then the $\varphi_{k-1}$-coordinate of $\chi_{k}(e_1)$ belongs to $\{(B_\ell,+),(B_\ell,-),(I,0)\}$. However, the $\varphi_{k-1}$-coordinate of $\chi_{k}(e_2)$ belongs to $\{(B_{\ell-1},+),(B_{\ell-1},-)\}$. A contradiction.
If $\ell=2$, then the $\varphi_{k-1}$-coordinate of $\chi_{k}(e_1)$ belongs to $\{(B_2,+),(B_2,-)\}$. However, the $\varphi_{k-1}$-coordinate of $\chi_{k}(e_2)$ belongs to $\{(D,0),(A_\ell,0)\}$ for some $\ell\in [2,k-2]$. Again a contradiction.

This completes the induction step and so the assertion follows. \hfill$\Box$

\medskip

Now, Theorem \ref{center-2} follows from the following lemma.
\begin{lemma}\label{(k+2,3)-coloring}
Let $\chi_k$ ($k\ge4$) be the edge-coloring of $H_k$ as above. Then $\chi_k$ is a $(k+2,3)$-coloring.
\end{lemma}
\textbf{Proof.} The proof is by using the induction on $k\ge4$. The base case where $k=4$ holds from Lemma \ref{(6,3)-coloring}. In general, for $k>4$, suppose the assertion holds for $k-1$ and we will show it also holds for $k$.

Let $X_k=\{a_1,\ldots ,a_{k+2}\}$ be vertices of $H_k$ forming a $K_{k+2}^{(k)}$ with $a_1<\cdots <a_{k+2}$, we will show that there has at least three colors in edges of $X_k$. Recall that $\delta_i=\delta(a_i,a_{i+1})$ for $i\in[ k+1]$, and $\delta_{ij}=\delta(a_i,a_j)$ for $1\le i<j\le k+2$. Let $e_1=(a_1,\ldots,a_k)$, $e_2=(a_2,\ldots,a_{k+1})$, and $e_3=(a_3,\ldots,a_{k+2})$ be the three edges in $X_k$, then $\delta(e_1)=\{\delta_1,\ldots,\delta_{k-1}\}$, $\delta(e_2)=\{\delta_2,\ldots,\delta_{k}\}$, and $\delta(e_3)=\{\delta_3,\ldots,\delta_{k+1}\}$.

Suppose that $\chi_k(e_1)$, $\chi_k(e_2)$, and $\chi_k(e_3)$ are distinct to each other, then this gives three colors to the edges in $X_k$ and we need do nothing. Therefore, there are at least two of them are equal. We split the proof into three cases as follows.

\medskip
\emph{Case 1:} $\chi_k(e_1)=\chi_k(e_2)$.

\medskip
For this case, $\varphi_{k-1}(\delta_1,\ldots,\delta_{k-1})=\varphi_{k-1}(\delta_2,\ldots,\delta_k)$.
If $\delta(X_k)$ is monotone, then the $K_{k+1}^{(k-1)}$ on $\{\delta_1,\ldots,\delta_{k+1}\}$ has three colors since $\chi_{k-1}$ is a $(k+1,3)$-coloring from the induction hypothesis and this gives at least three colors to the edges in $X_k$. Thus, $\delta(X_k)$ has a local extremum. Let $\delta_\ell$ be the first local extremum, where $\ell\in [2,k]$.

Suppose $\ell\in [3,k-1]$. We may assume $\delta_\ell$ is a local maximum without loss of generality, then
$\varphi_{k-1}(\delta_1,\ldots,\delta_{k-1})\in\{(B_\ell,+),(B_\ell,-),(I,0)\}$, and $\varphi_{k-1}(\delta_2,\ldots,\delta_k)\in\{(B_{\ell-1},+),(B_{\ell-1},-)\}$. This leads to a contradiction.
If $\ell=2$, then we may assume that $\delta_2$ is a local maximum without loss of generality, i.e., $\delta_1<\delta_2>\delta_3\cdots$. Then,
$\varphi_{k-1}(\delta_1,\ldots,\delta_{k-1})\in\{(B_2,+),(B_2,-)\}$, and $\varphi_{k-1}(\delta_2,\ldots,\delta_k)\in\{(D,0),(A_{j},0)\}$ for some $j \in [2, k-2]$.
A contraction. If $\ell=k$, then we may assume that $\delta_k$ is a local maximum without loss of generality, i.e., $\delta_1<\cdots<\delta_{k-1}<\delta_k>\delta_{k+1}$. Then the $K_k^{(k-1)}$ on $\{\delta_1,\ldots,\delta_k\}$ has two colors since $\chi_{k-1}$ is a $(k+1,2)$-coloring from Lemma \ref{(k+1,2)-coloring} and this gives two colors to the edges of $X_k$ within $\{a_1,\ldots,a_{k+1}\}$ and the $\varphi_{k-1}$-coordinates are $(I,0)$. Moreover, $\varphi_{k-1}(\delta_3,\ldots,\delta_{k+1})\in\{ (B_{k-2},+),(B_{k-2},-)\}$, so there are at least three colors on $X_k$.

\medskip
\emph{Case 2:} $\chi_k(e_2)=\chi_k(e_3).$

\medskip
If $\delta(X_k)$ is monotone, then we are done similarly as in Case 1. Thus, let $\delta_\ell$ be the first local extremum, where $\ell\in [2,k]$. If $\ell\in[3,k]$, then we are also done as in Case 1. So we may assume $\ell=2$ and $\delta_2$ is a local maximum without loss of generality, i.e., $\delta_1<\delta_2>\delta_3\cdots$. Moreover, $\varphi_{k-1}(\delta_2,\ldots,\delta_k)=\varphi_{k-1}(\delta_3,\ldots,\delta_{k+1})$ from $\chi_k(e_2)=\chi_k(e_3)$, then we have $\delta_2>\delta_3>\cdots>\delta_{k+1}$. By a similar argument as above, there are at least three colors on $X_k$ by using Lemma \ref{(k+1,2)-coloring}.

\medskip
\emph{Case 3:} $\chi_k(e_1)=\chi_k(e_3).$

\medskip

Similarly as above, if $\delta(X_k)$ is monotone then we are done. Thus, let $\delta_\ell$ be the first local extremum, where $\ell\in [2,k]$. If $\ell\in[3,k]$, then we are also done similarly as above. So we may assume $\ell=2$, and $\delta_2$ is a local maximum without loss of generality, i.e., $\delta_1<\delta_2>\delta_3\cdots$. Moreover, $\varphi_{k-1}(\delta_1,\ldots,\delta_{k-1})=\varphi_{k-1}(\delta_3,\ldots,\delta_{k+1})$ since $\chi_k(e_1)=\chi_k(e_3)$, and thus we have $\delta_3<\delta_4>\delta_5$. It follows from Property D that $\delta_1\neq \delta_3$ and $\delta_3\neq \delta_5$.

Suppose first $\delta_1<\delta_3$, i.e., $\delta_1<\delta_3<\delta_4>\delta_5\cdots$. Consider $e_1'=(a_1,a_3,\ldots,a_{k+1})$, then $\delta(e_1')=\{\delta_1,\delta_3,\ldots,\delta_k\}$ since $\delta_{13}=\delta_1$ from Property C. Thus, we have $\varphi_{k-1}(\delta_1,\delta_3,\ldots,\delta_k)\in \{(B_3,+),(B_3,-),(I,0)\}$. Recall $\varphi_{k-1}(\delta_1,,\ldots,\delta_{k-1})=(B_2,+)$ and $\varphi_{k-1}(\delta_2,\ldots,\delta_k)=(A_2,0)$. Thus, the $\varphi_{k-1}$-coordinates of $\chi_k(e_1')$, $\chi_k(e_1)$ and $\chi_k(e_2)$ are distinct, and so we have at least three colors on $X_k$. Thus, we may assume $\delta_1>\delta_3$.

Suppose $\delta_3>\delta_5$, i.e., $\delta_2>\delta_3>\delta_5\cdots$. Let $e_1'=(a_2,a_3,a_5\ldots,a_{k+2})$, and we have that $\delta(e_1')=\{\delta_2,\delta_3,\delta_5,\ldots,\delta_{k+1}\}$ since $\delta_{35}=\delta_3$ from Property C. Then, $\varphi_{k-1}(\delta_2,\delta_3,\delta_5,\ldots,\delta_{k+1})\in\{(D,0),(A_j,0)\}$ for some $j\in[3,k-2]$. Recall $\varphi_{k-1}(\delta_2,\ldots,\delta_k)=(A_2,0)$ and $\varphi_{k-1}(\delta_3,\ldots,\delta_{k+1})=(B_2,-)$. Thus, the $\varphi_{k-1}$-coordinates of $\chi_k(e_1')$, $\chi_k(e_2)$ and $\chi_k(e_3)$ are distinct, and so we have at least three colors on $X_k$. Thus, we may assume $\delta_3<\delta_5$.

From the above, we conclude that $\delta_1>\delta_3<\delta_5$. Therefore, $\varphi_{k-1}(\delta_1,\ldots,\delta_{k-1})=(B_2,-)$ and $\varphi_{k-1}(\delta_3,\ldots,\delta_{k+1})=(B_2,+)$. Recall $\varphi_{k-1}(\delta_2,\delta_3,\ldots,\delta_k)=(A_2,0)$. It follows that the $\varphi_{k-1}$-coordinates of $\chi_k(e_1)$, $\chi_k(e_2)$ and $\chi_k(e_3)$ are distinct, and thus there are again at least three colors on $X_k$.

This completes the proof.\hfill$\Box$

\bigskip

\noindent\emph{Remark.}~It is worth noting that the settings for $(B_\ell,+)$ and $(B_\ell,-)$ are only to guarantee that the case where $\delta_2$ is a local maximum in Case 3 can proceed. We don't need to classify $(A_\ell,0)$ like $\{(B_\ell,+),(B_\ell,-)\}$, indeed, if $\delta_2$ is a local minimum, then $\delta_1>\delta_2<\delta_3>\delta_4<\delta_5\cdots$. If $\delta_2<\delta_4$, then set $e_1'=(a_2,a_4,\ldots,a_{k+2})$. If $\delta_2>\delta_4$, then set $e_1'=(a_1,a_2,a_4,\ldots,a_{k+1})$. In this way, we can prove that $X_k$ has at least three colors by using similar arguments as above.

\end{spacing}
\end{document}